\documentstyle{amsppt}
\topmatter
\loadmsbm
\UseAMSsymbols
\hoffset=0.75truein
\voffset=0.5truein

\def\=#1{\accent"16 #1}

\def\refz{\relax}

\def\pt{\hbox{\it pt}}

\def\={\relax}

\def\i{I}
\def\was#1{\relax}

\def\diag|#1|#2|{\vbox to #1in {\vskip.3in\centerline{\tt Diagram #2}\vss} }
\def\v{\hskip -3.5pt }
\def\gram|#1|#2|#3|{
        {
        \smallskip
        \hbox to \hsize
        {\hfill
        \vrule \vbox{ \hrule \vskip 6pt \centerline{\it Diagram #2}
         \vskip #1in %
             \includegraphics{#3}\hrule }
        \v\vrule\hfill
        }
\smallskip}}

\title An Overview of the Kepler Conjecture\endtitle
\author Thomas C. Hales\endauthor
\endtopmatter

\document
\footnote""{\hfill version -- 7/29/98, revised 12/5/01}
\footnote""{\hfill Research supported in part by the NSF}

\bigskip
\head \refz 1. Introduction and Review\endhead

The series of papers in this volume gives a proof of the Kepler
conjecture, which asserts that the density of a packing of congruent
spheres in three dimensions is never greater than $\pi/\sqrt{18}\approx
0.74048\ldots$.
This is the oldest problem in discrete geometry
and is an important part of Hilbert's 18th problem.
An example of a packing achieving this density is
the face-centered cubic packing.

A packing of spheres is an arrangement of
nonoverlapping spheres of radius 1 in Euclidean space.
Each sphere is determined by its center, so equivalently it is a collection
of points in Euclidean space separated by distances of at least 2.
The density of a packing is defined as the $\limsup$ of
the densities of the partial packings formed by spheres inside
a ball with fixed center of radius $R$.
(By taking the $\limsup$,
rather than $\liminf$ as the density, we prove the Kepler
conjecture in the strongest possible sense.)
Defined as a limit, the density
is insensitive to changes in the packing in any bounded region.
For example, a finite number of spheres can be removed from the
face-centered cubic packing without affecting its density.

Consequently, it is not possible to hope for any strong uniqueness results
for packings of optimal density.  The uniqueness established by this
work is as strong as can be hoped for.
It shows that certain local
structures (decomposition stars)
attached to the face-centered cubic (fcc) and hexagonal-close
packings (hcp)
are the only structures that maximize a local density function.

Although we do not pursue this point, Conway and Sloane develop
a theory of tight packings that is more restrictive than having the
greatest possible density \cite{CoSl95}.
An open problem is to prove that their
list of tight packings in three dimensions is complete.

The face-centered cubic packing appears in
Diagram 1.1.

\smallskip
\gram|2.2|1.1|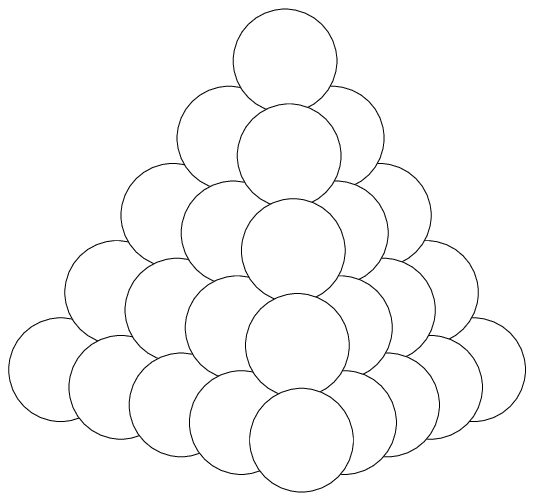|
\smallskip

\head \refz 2. Early History, Hariot, and Kepler\endhead

The study of the mathematical properties of the face-centered cubic
packing can be traced back to a Sanskrit work composed around 499 CE.
I quote an extensive passage from the commentary that K. Plofker
has made about the formula
for the number of balls in triangular piles\cite{Plo00}:
\medskip

{
\narrower
\font\ninerm=cmr9
\ninerm
\def\=#1{\accent"16 #1}

 The excerpt below is taken
 from a Sanskrit work composed around 499 CE, the \=Aryabha\d t\={\i}ya
 of \=Aryabha\d ta, and the commentary on it written in 629 CE
 by Bh\=askara (I).  The work is a compendium of various rules in
 mathematics and mathematical astronomy, and the results are probably
 not due the \=Aryabha\d ta himself but derived from an earlier source:
 however, this is the oldest source extant for them.  (My translation's
 from the edition by K. S. Shukla, {\it The \=Aryabha\d t\={\i}ya of
 \=Aryabha\d ta with the Commentary of Bh\=askara I and Some\'svara},
 New Delhi: Indian National Science Academy 1976; my inclusions are
 in square brackets. There is a corresponding English translation
 by Shukla and K. V. Sarma, {\it The \=Aryabha\d t\={\i}ya
 of \=Aryabha\d ta}, New Delhi: Indian National Science Academy 1976.
 It might be easier to get hold of the earlier English translation by
 W. E. Clark, {\it The \=Aryabha\d t\={\i}ya of \=Aryabha\d ta},
 Chicago: University of Chicago Press, 1930.)

      Basically, the rule considers the series in arithmetic
 progression
 $
 S_i = 1 + 2 + 3 + \ldots + i
 $
 (for whose sum the formula is known) as the number of objects
 in the $i$th layer of a pile with a total of $n$ layers, and
 specifies the following two equivalent formulas for the ``accumulation
 of the pile'' or $ \sum_{i=1}^n S_i $:

 $$
 \sum_{i=1}^n S_i = {{n(n+1)(n+2)}\over 6},
 $$

 $$
 \sum_{i=1}^n S_i = {{(n+1)^3 - (n+1)}\over 6}.
 $$

 What he says is this:

 {\it \=Aryabha\d t\={\i}ya}, Ga\d nitap\=ada 21:

 {\narrower
    For a series [lit. ``heap''] with a common difference and first term
 of 1, the product of three [terms successively] increased by 1 from
 the total, or else the cube of [the total] plus 1 diminished by [its]
 root, divided by 6, is the total of the pile [lit. ``solid heap''].
 }

 Bh\=askara's commentary on this verse:

 {\narrower
    [This] heap [or] series is specified as having one for its common
 difference and initial term. This same series with one for its common
 difference and initial term is said [to be] ``heaped up.'' ``The
 product of three [terms successively] increased by one from the total''
 of this so-called heaped-up ``series with one for its common
 difference and initial term'': i.e., the product of three terms, starting
 from the total and increasing by one. Namely, the total, that plus one,
 and [that] plus one again. That [can] be stated [as follows]: the total,
 that plus one, and that total plus two. The product of those three
 divided by 6 is the ``solid heap,'' the accumulation of the series.
 Now another method: The cube of the root equal to that [total] plus
 one is diminished by its root, and divided by 6: thus it follows.
 ``Or else'': [i.e.], the cube of that root plus one, diminished by
 its own root, divided by 6, is the ``solid heap.''
    Example: Series with 5, 8, and 14 respectively for their total layers:
 tell me [their] triangular-shaped piles.
    In order, the totals are 5, 8, 14.
    Procedure: Total 5. This plus one: 6. This plus one again: 7. Product
 of those three: 210. This divided by 6 is the accumulation of the series:
 35.
    [He goes on to give the answers for the second two cases, but you
 doubtless get the picture.]   --  K. Plofker
 }

}

The modern mathematical study of spheres and their close packings can be
traced to T. Hariot.  Hariot's work -- unpublished, unedited,
and largely undated -- shows a preoccupation with sphere packings.
He seems to have first taken an interest in packings at the
prompting of Sir Walter Raleigh.  At the time, Hariot was Raleigh's
mathematical assistant,  and Raleigh gave him the problem of determining
formulas for the number of cannonballs in regularly stacked piles.
In 1591 he prepared a chart of triangular numbers for Raleigh.
Shirley, Hariot's biographer, writes,

{
\narrower
\font\ninerm=cmr9
\ninerm

    Obviously, this is a quick reference chart prepared for Ralegh to give information on the ground space required for the storage of cannon
    balls in connection with the stacking of armaments for his marauding vessels. The chart is ingeniously arranged so that it is possible to
    read directly the number of cannon balls on the ground or in a pyramid pile with triangular, square, or oblong base. All of this Harriot had
    worked out by the laws of mathematical progression (not as Miss Rukeyser suggests by experiment), as the rough calculations
    accompanying the chart make clear. It is interesting to note that on adjacent sheets, Harriot moved, as a mathematician naturally would,
    into the theory of the sums of the squares, and attempted to determine graphically all the possible configurations that discrete particles
    could assume -- a study which led him inevitably to the corpuscular or atomic theory of matter originally deriving from Lucretius and
    Epicurus. \cite{Shi83,p.242}

}

\smallskip
Hariot connected sphere packings to Pascal's triangle long before Pascal
introduced the triangle.
See Diagram 2.1.

\smallskip
\gram|5.0|2.1 Hariot's triangle|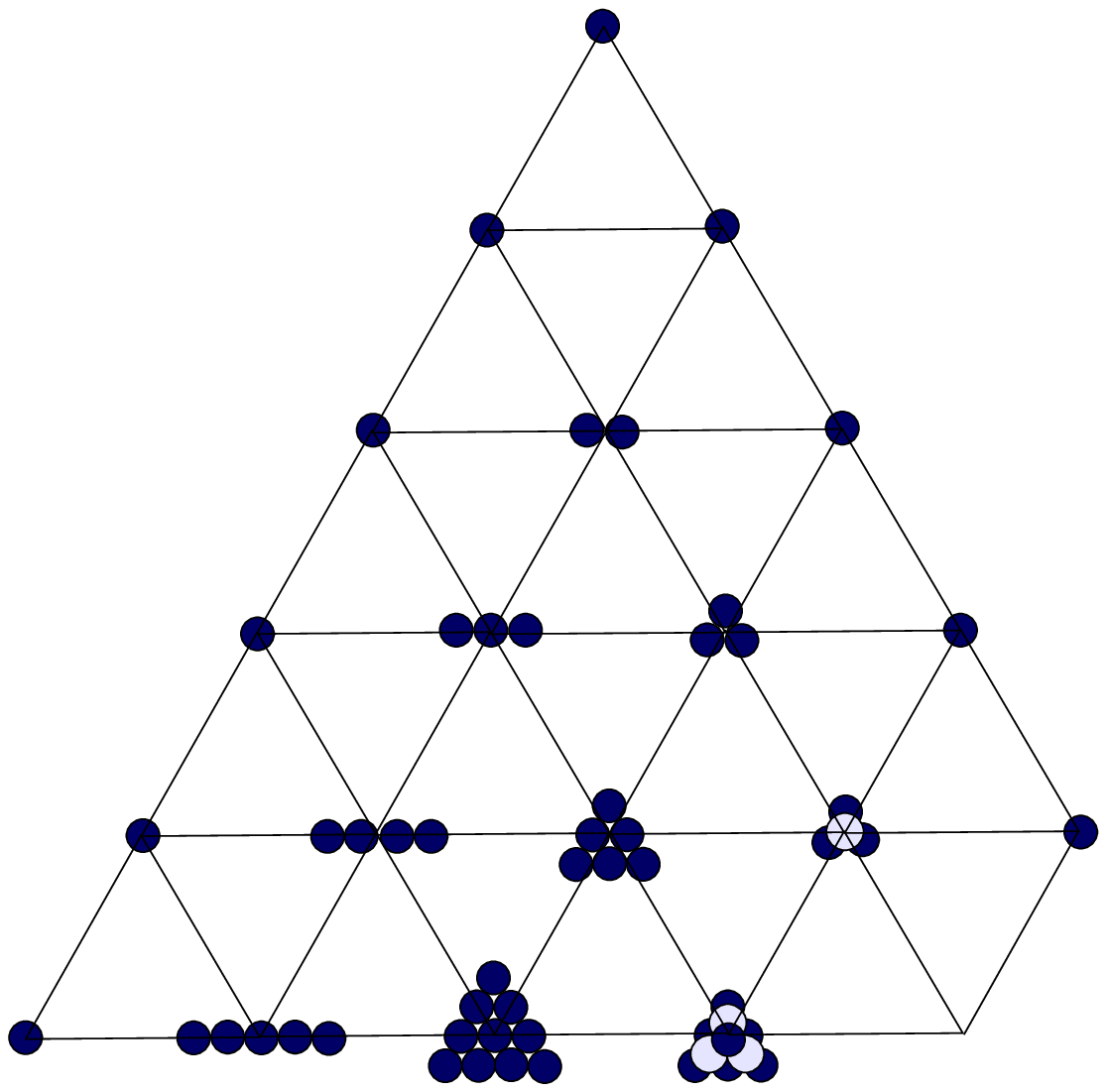|
\smallskip

Hariot was the first to distinguish between the
face-centered cubic and hexagonal close packings \cite{Mas66,p.52}.

Kepler became involved in sphere packings through his correspondence
with Hariot in the early years of the 17th century.
Kargon writes, in his history of atomism in England,

{
\narrower
\font\ninerm=cmr9
\ninerm

    Hariot's theory of matter appears to have been virtually that of Democritus, Hero of Alexandria, and, in a large measure, that of Epicurus
    and Lucretius. According to Hariot the universe is composed of atoms with void space interposed. The atoms themselves are eternal and
    continuous. Physical properties result from the magnitude, shape, and motion of these atoms, or corpuscles compounded from them$\ldots$.

    Probably the most interesting application of Hariot's atomic theory was in the field of optics. In a letter to Kepler on 2 December 1606
    Hariot outlined his views. Why, he asked, when a light ray falls upon the surface of a transparent medium, is it partially reflected and
    partially refracted? Since by the principle of uniformity, a single point cannot both reflect and transmit light, the answer must lie in the
    supposition that the ray is resisted by some points and not others.

    ``A dense diaphanous body, therefore, which to the sense appears to be continuous in all parts, is not actually continuous. But it has
    corporeal parts which resist the rays, and incorporeal parts vacua which the rays penetrate$\ldots$''

    It was here that Hariot advised Kepler to abstract himself mathematically into an atom in order to enter `Nature's house'. In his reply of 2
    August 1607, Kepler declined to follow Harriot, ad atomos et vacua. Kepler preferred to think of the reflection-refraction problem in terms
    of the union of two opposing qualities --
    transparence and opacity. Hariot was surprised. ``If those assumptions and reasons satisfy you, I
    am amazed.'' \cite{Kar66,p.26}

}

\smallskip
Despite Kepler's initial reluctance to adopt an atomic theory, he
was eventually swayed, and in 1611 he published an essay that explores
the consequences of a theory of matter composed of small
spherical particles.  Kepler's essay was the
``first recorded step towards a mathematical theory of the genesis
of inorganic or organic form'' \cite{Why66,p.v}.

Kepler's essay describes
the face-centered cubic packing and asserts that ``the packing will
be the tightest possible, so that in no other arrangement  could more
pellets be stuffed into the same container.''  This assertion has
come to be known as the Kepler conjecture.   The purpose of this
collection of papers is to give a proof of this conjecture.

\head 3. History \endhead

The next episode in the history of this problem is a debate between
Isaac Newton and David Gregory.
  Newton and Gregory discussed the question of how many spheres
of equal radius
can be arranged to touch a given sphere.  This is the three-dimensional
analogue of the simple fact that in two dimensions six pennies,
but no more, can
be arranged to touch a central penny.  This is the kissing-number
problem in $n$-dimensions.
In three dimensions, Newton said that the
maximum was 12 spheres, but Gregory claimed that 13 might be possible.

Newton was correct.
In the 19th century, the first papers claiming a proof of the
kissing-number problem appeared
in \cite{Ben74}, \cite{Gun75}, \cite{Hop74}.
Although some writers cite these papers
as a proof, they are hardly rigorous by today's standards.
Another incorrect
proof appears in \cite{Boe52}.
  The first proper proof was obtained
by B. L. van der Waerden and Sch\"utte in 1953 \cite{Sch53}.
An elementary proof appears in Leech \cite{Lee56}.
The influence of van der Waerden, Sch\"utte, and Leech upon the
papers in this collection is readily apparent.  Although the
connection between the Newton-Gregory problem and Kepler's problem
is not obvious, L. Fejes T\'oth in 1953, in the first
work describing a strategy to prove the Kepler conjecture, made
a quantitative version of the Gregory-Newton problem the first step
\cite{Fej53}.

The two-dimensional analogue of the Kepler conjecture is to show
that the honeycomb packing in two dimensions gives the highest
density.  This result was established in 1892 by Thue, with
a second proof appearing in 1910 (\cite{Thu92}, \cite{Thu10}).
A number of other proofs have appeared since then.  Three are particularly
notable.
Rogers's proof
generalizes to give a bound on the density of
packings in any dimension \cite{Rog58}.
A proof by L. Fejes T\'oth extends to give
bounds on the density of packings of convex disks \cite{Fej50}.
A third proof, also by L. Fejes T\'oth, extends to
non-Euclidean geometries \cite{Fej53}.
Another early proof appears in \cite{SeM44}.

In 1900, Hilbert made the Kepler conjecture part of his 18th problem
\cite{Hil01}.
Milnor, in his review of Hilbert's 18th problem, breaks the problem
into three parts \cite{Mil76}.

{
\narrower
\font\ninerm=cmr9
\ninerm

1.  Is there in $n$-dimensional Euclidean Space $\ldots$ only a finite
number of essentially different kinds of groups of motions with a
[compact] fundamental region?

2.  Whether polyhedra also exist which do not appear as fundamental
    regions of groups of motions, by means of which nevertheless
    by a suitable juxtaposition of congruent copies a complete filling
    up of all [Euclidean] space is possible?

3.  How can one arrange most densely in space an infinite number
    of equal solids of given form, e.g. spheres with given radii $\ldots$,
    that is, how can one so fit them together that the ratio of the
    filled to the unfilled space may be as great as possible?

}

\smallskip
Writing of the third part, Milnor states,

{
\narrower
\font\ninerm=cmr9
\ninerm

For $2$-dimensional disks
this problem has been solved by Thue and Fejes T\'oth, who showed
that the expected hexagonal (or honeycomb) packing of circular disks
in the plane is the densest possible.  However, the corresponding
problem in $3$ dimensions remains unsolved.  This is a scandalous
situation since the (presumably) correct answer has been known since
the time of Gauss. (Compare Hilbert and Cohn-Vossen.)  All that is
missing is a proof.

}

\head 4. The Literature \endhead

Past progress toward the Kepler conjecture can be arranged into
four categories:
(1) bounds on the density, (2) descriptions of classes of packings
for which the bound of $\pi/\sqrt{18}$ is known,
(3) convex bodies
other than spheres for which the packing density can be determined
precisely, (4) strategies of proof.

\subhead 4.1. Bounds\endsubhead

Various upper bounds have been established on the density of
packings.
\smallskip

{\obeylines

0.884 \cite{Bli19},
0.835 \cite{Bli29},
0.828 \cite{Ran47},
0.7797 \cite{Rog58},
0.77844 \cite{Lin86},
0.77836 \cite{Mud88},
0.7731 \cite{Mud93}.

}

\smallskip
Rogers's is a particularly natural bound.
  As the dates indicate, it remained the best available
bound for many years.  His monotonicity lemma and his
decomposition of Voronoi cells into simplices have become important
elements in the proof of the Kepler conjecture.
We give a new proof of Rogers's bound
in ``Sphere Packings III.''  A function $\tau$,
used throughout this
collection, measures the departure of various objects from
Rogers's bound.

Muder's bounds, although they appear to be rather small
improvements of Rogers's bound, are the first to make use
of the full Voronoi cell in the the determination of densities.
As such, they mark a transition to a greater level of
sophistication and difficulty.  Muder's influence on the work
in this collection is also apparent.

A sphere packing admits a Voronoi decomposition: around
every sphere take the convex region consisting of points closer to that sphere
center than to any other sphere center.   L. Fejes T\'oth's
dodecahedral
conjecture asserts that the Voronoi cell of smallest volume is
a regular dodecahedron with inradius 1 \cite{Fej42}.
The dodecahedral conjecture implies a bound of 0.755 on sphere
packings.  L. Fejes T\'oth actually gave a complete proof except
for one estimate. A footnote in his paper documents the gap, ``In the
proof, we have relied to some extent solely on intuitive
observation [Anschauung].''
 As L. Fejes T\'oth pointed out, that estimate is extraordinarily
difficult, and the dodecahedral conjecture has resisted all efforts
until now \cite{McL98}.

The missing estimate in L. Fejes T\'oth's paper is an explicit form
of the Newton-Gregory problem.  What is needed is an explicit bound
on how close the 13th sphere can come to touching the central
sphere.  Or more generally, minimize the sum of the distances
of the 13 spheres from the central sphere.
No satisfactory bounds are known.  Boerdijk has a conjecture for the arrangement
that minimizes the average distance of the 13 spheres from the
central sphere.   Van der Waerden
has a conjecture for the closest arrangement of 13 spheres in which
all spheres have the same distance from the central sphere.
Bezdek has shown that the dodecahedral conjecture would follow from
weaker bounds than those originally proposed by L. Fejes T\'oth
\cite{Bez97}.

A proof of the dodecahedral conjecture has traditionally been
viewed as the first step toward a proof of the Kepler conjecture,
and if little progress has been made until now toward a complete
solution of the Kepler conjecture, the difficulty of the dodecahedral
conjecture is certainly responsible to a large degree.

\subhead 4.2. Classes of packings\endsubhead

If the infinite dimensional space of all packings is too unwieldy,
we can ask if it is possible to establish the bound $\pi/\sqrt{18}$
for packings with special structures.

If we restrict the problem
to packings whose sphere centers are the points of a lattice, the
 packings are described by a finite number of parameters, and the
problem becomes much more
accessible.  Lagrange proved that the densest lattice packing
in two dimensions is the familiar honeycomb arrangement \cite{Lag73}.
Gauss proved that
the densest lattice packing in three dimensions is the face-centered
cubic \cite{Gau31}.
In dimensions 4--8, the optimal lattices
are described by their root systems,
$A_2$, $A_3$, $D_4$, $D_5$, $E_6$, $E_7$, and $E_8$.
A. Korkine and G. Zolotareff showed that $D_4$
and $D_5$ are the densest lattice packings
in dimensions 4 and 5 (\cite{KoZ73}, \cite{KoZ77}).
Blichfeldt determined the densest lattice packings in
dimensions 6--8 \cite{Bli35}.
Beyond
dimension $8$, there are no proofs of optimality, and yet there
are many excellent candidates for the densest lattice packings
such as the Leech lattice in dimension $24$.
For a proof of the existence of optimal lattices, see \cite{Oes90}.

Although lattice packings are of particular interest because they
relate to so many different branches of mathematics, Rogers has
conjectured that in sufficiently high dimensions, the densest
packings are not lattice packings \cite{Rog64}.   In fact,
the densest known packings in various dimensions are not lattice
packings.  The third edition of \cite{CoSl93} gives several examples
of nonlattice packings that are denser than any known lattice
packings (dimensions 10, 11, 13, 18, 20, 22).
The densest packings of typical convex sets in the plane,
in the sense
of Baire categories, are not lattice packings \cite{Fej95}.

Gauss's theorem on lattice densities has been generalized by
A. Bezdek, W. Kuperberg, and E. Makai, Jr. \cite{BKM91}.
They showed that packings of parallel
strings of spheres never have density greater than $\pi/\sqrt{18}$.

\subhead 4.3. Other convex bodies\endsubhead

If the optimal sphere packings are too difficult to determine,
we might ask whether
the problem can be solved for other convex bodies.
To avoid trivialities, we restrict our attention to convex bodies
whose packing density is strictly less than 1.

  The first convex body in Euclidean 3-space that does not tile
for which the packing density was explicitly determined is
an infinite cylinder \cite{Bez90}.
Here A. Bezdek and W. Kuperberg prove
that the
optimal density is obtained by arranging the cylinders in
parallel columns in the honeycomb arrangement.

In 1993, J. Pach exposed the humbling depth of our ignorance when he issued
the challenge to determine the packing density for some bounded convex
body that does not tile space \cite{MP93}.
(Pach's question is more revealing than anything I can write on
the subject of discrete geometry.)
 This question was answered by
A. Bezdek \cite{Bez94}, who determined the packing density of a rhombic
dodecahedron that has one corner clipped so that it no longer tiles.
The packing density equals the ratio of the
volume of the clipped
rhombic dodecahedron to the volume of the unclipped rhombic dodecahedron.

\subhead 4.4. Strategies of proof\endsubhead

In 1953, L. Fejes T\'oth proposed a program to prove the
Kepler conjecture \cite{Fej53}.
A single Voronoi cell cannot lead to a bound better
than the dodecahedral conjecture.   L. Fejes T\'oth considered
weighted averages of the volumes of collections of Voronoi cells.
 These weighted
averages involve up to 13 Voronoi cells.  He showed that if a particular
weighted average of volumes is greater than the volume of the
rhombic dodecahedron, then the Kepler conjecture follows.
The Kepler conjecture is an optimization problem in an infinite
number of variables.  L. Fejes T\'oth's weighted-average argument
was the first indication that it might be possible to reduce
the Kepler conjecture to a problem in a finite number of variables.
Needless to say, calculations involving the weighted averages of the
volumes of
several Voronoi cells will be significantly more difficult than those
involved in establishing the dodecahedral conjecture.

To justify his approach, which limits the number of Voronoi cells
to 13, Fejes T\'oth needs a preliminary estimate of how close
a 13th sphere can come to a central sphere.  It is at this point
in his formulation of the Kepler conjecture that an explicit
version of the Newton-Gregory problem is required.  How
close can 13 spheres come to a central sphere, as measured by
the sum of their distances from the central sphere?

Strictly speaking, neither L. Fejes T\'oth's program nor my own program
reduces the Kepler
conjecture to a finite number of variables, because if
it turned out that one of
the optimization problems in finitely many
variables had an unexpected global maximum, the program would
fail, but the Kepler conjecture would remain intact.
In fact, the failure of a program has
no implications for the Kepler conjecture.  The proof that the
Kepler conjecture reduces to a finite number of variables comes only as
corollary to the full proof of the Kepler conjecture.

L. Fejes T\'oth made another significant suggestion in \cite{Fej64}.
He was the first to suggest the use of computers in the Kepler conjecture.
After describing his program, he writes,

{
\narrower
\font\ninerm=cmr9
\ninerm

Thus it seems that the problem can be reduced to the determination
of the minimum of a function of a finite number of variables,
providing a programme realizable in principle.  In view of the
intricacy of this function we are far from attempting to
determine the exact minimum.  But, mindful of the rapid development
of our computers, it is imaginable that the minimum may
be approximated with great exactitude.

}

\smallskip
The most widely publicized attempt to prove the Kepler conjecture
was that of Wu-Yi Hsiang \cite{Hsi93}.  (See also
\cite{Hsi93a}, \cite{Hsi93b}.)  Hsiang's approach can be
viewed as a continuation and extension of L. Fejes T\'oth's
program.  Hsiang's paper contains major
gaps and errors \cite{CoHMS94}.
  The mathematical arguments against his argument appear
in my
debate with him in the {\it Mathematical Intelligencer}
(\cite{Hal94}, \cite{Hsi95}).
There are now many published sources that agree with the central
claims of \cite{Hal94} against Hsiang.
Conway and Sloane report that the paper ``contains serious flaws.''
G. Fejes T\'oth feels that ``the greater part of the work has yet
to be done'' \cite{Fej95}.   K. Bezdek concluded,
after an extensive study of Hsiang's work, ``his work is far from being
complete and correct in all details'' \cite{Bez97}.
 D. Muder writes, ``the
community has reached a consensus on it: no one buys it'' \cite{Mud97}.

\head 5. Experiments with other Decompositions\endhead

  The results of my early efforts to prove the Kepler conjecture
are published in \cite{Hal92}, \cite{Hal93}.
A decomposition dual to the Voronoi decomposition is the
Delaunay decomposition.   The papers show that the Delaunay decomposition
also leads to an optimization problem in a finite number of variables
that implies the Kepler conjecture.   This approach led to the first
proof that the face-centered cubic and hexagonal-close packings are
locally optimal in the appropriate sense.  Also, the optimization problem
was studied numerically, and this gave the first direct numerical
evidence of the truth of the Kepler conjecture.  This method
also shows how the Delaunay and Voronoi decompositions can be superimposed
in a way that gives slightly better bounds than either does individually.

Unfortunately, the approach based on Delaunay decomposition also
rapidly runs into enormously difficult technical complications.
There was little hope of solving the problem using rigorous methods.
I set my program aside until 1994, when I began again with renewed
energy.   I proposed a new decomposition that is a hybrid
of the Voronoi and Delauanay decompositions.  Since then,
this hybrid decomposition
has passed through a long series of refinements.  The
most important stages of that development are described in
``Sphere Packings I'' and ``Sphere Packings II.''
Its final form was worked out in collaboration with
S. Ferguson in ``A Formulation of the Kepler conjecture.''

\head 6. Complexity\endhead

Why is this a difficult problem?  There are many ways to answer this
question.

This is an optimization problem in an infinite number
of variables.  In many respects, the central problem has been to
formulate a good finite dimensional approximation to
the density of a packing.  Beyond this, there remains an extremely
difficult problem in global optimization, involving
nearly 150 variables.
We recall that even very simple classes
of nonlinear optimization problems, such as quadratic optimization
problems,
are NP-hard \cite{HoPT95}.  A general highly nonlinear program of this
size is regarded by most researchers as hopeless (at least as far as rigorous
methods are concerned).

There is a considerable literature on many closely related nonlinear
optimization problems (the Tammes problem, circle packings, covering
problems, the Leonard-Jones potential, Coulombic energy minimization
of point particles, and so forth).
Many of our expectations about nonlattice packings are formed
by the extensive experimental data that have
been published on these problems.
The literature leads one to expect a rich abundance of critical
points, and yet it leaves one with a certain skepticism about the
possibility of establishing general results rigorously.

The extensive survey of circle packings
in \cite{Mel97} gives a broad overview of the progress and limits
of the subject.
Problems involving a few circles can be
trivial to solve.
Problems involving several circles in the plane can be solved with
sufficient ingenuity.
With the aid of computers, various  problems involving
a few more circles can be
treated by rigorous methods.
Beyond that, numerical methods
give approximations but no rigorous solutions.
Melissen's account of
the 20-year quest for the best separated arrangement of 10
points in a unit square is particularly revealing of the complexities
of the subject.

Kepler's problem has a particularly rich collection of
(numerical) local maxima that come uncomfortably close to the
global maximum \cite{Hal92}.
These local maxima
explain in part why a large number (around 5000) of planar maps
are generated as part of the proof of the conjecture.  Each
planar map leads to a separate nonlinear optimization problem.

\head 7. Contents of Papers\endhead

There are two papers that form part of the proof of the Kepler conjecture
that precede this collection.

``Sphere Packings I''
 defines a decomposition of space according
to a hybrid of the Delauany simplices and Voronoi decomposition.
We will not go into the computational difficulties that arise with the
pure Delaunay and pure Voronoi decompositions.  Suffice it to
say that the technical problems surrounding either of the pure
strategies are immense.

The parts of the hybrid decomposition
that lie around a given sphere center form a star-shaped region
called the decomposition star.
Just as a bound on the volume of
Voronoi cells leads to a bound on the density of a packing,
a bound on a function called the score on the space of decomposition
stars leads to a bound on the density.   A sufficiently good
bound on the score $8\,\pt\approx 0.4429$
leads to the desired bound $\pi/\sqrt{18}$
on the density of a packing.

The aim is to show that the bound $8\,\pt$ holds for every decomposition star.
With each decomposition star is associated a planar map that describes
its combinatorial structure.
``Sphere Packings I''
proves the bound $8\,\pt$ for every decomposition star
whose planar map is a triangulation.

``Sphere Packings I''
 also lays the foundation for other papers by
listing classical formulas for
dihedral angles, solid angles, volumes of Voronoi cells, and
so forth.  It uses interval arithmetic to prove various inequalities
by computer.

``Sphere Packings II''
 proves a local optimality result.  The decomposition
stars of the face-centered cubic and hexagonal close packings score
$8\,\pt$. If they are deformed, the score drops below $8\,\pt$.
``Sphere Packings II'' proves local optimality in the strong sense
that among the planar graphs that occur in the fcc packing and the
hcp, the only decomposition stars that score $8\,\pt$ are
those of the fcc and hcp.  That is, they are the only global maxima
on the connected components of the space of decomposition stars to
which they belong.

The contributions to this collection begin with the paper
``A Formulation of the Kepler Conjecture.'' It
brings the structure of the  decomposition star and the scoring functions
into final form.
The papers in the series were written over a five-year period.  As
our investigations progressed, we found that it was necessary to
make some adjustments.  However, we had no desire to start over,
abandoning the results of ``Sphere Packings I''
 and ``Sphere Packings II.''
``A Formulation'' gives  a new decomposition of space.  The Delaunay
simplices are replaced with simplices that
approximate them (quasi-regular tetrahedra
and quarters).  Voronoi cells are replaced as well with slightly
modified objects called $V$-cells.  The scoring function is also
adjusted.   ``A Formulation''
shows that all of the main theorems from ``Sphere Packings I'' and
``Sphere Packings II'' can be easily recovered in this new context with a few
simple lemmas.   It also lays the foundations for much of
``Sphere Packings III,'' ``Sphere Packings IV,'' and ``Sphere Packings V.''

``Sphere Packings V'' treats the decomposition stars that have a
particular planar map, the pentagonal prisms.   This arrangement
of spheres has a long history.  Boerdijk was the first to realize
its importance, and used it to produce a counterexample to a conjecture
of L. Fejes T\'oth about an explicit form of the Gregory-Newton
problem \cite{Boe52}.  This arrangement shows up in various other places, such
as \cite{SHDC95}
and \cite{Hsi93}.  This arrangement gives a counterexample to
my earliest expectations, that a pure Delauany decomposition should
lead directly to a finite dimensional formulation of the Kepler conjecture.
This arrangement is either a counterexample or comes extremely close
to being a counterexample in the various formulations of the Kepler
conjecture that we have considered.  It had to be singled out for
special and careful attention.  The purpose of
``Sphere Packings V''
is to show that this particular type of decomposition star has score
under $8\,\pt$.

``Sphere Packings III,'' ``Sphere Packings IV,'' and a final paper,
``The Kepler Conjecture,'' treat all the remaining possibilities.
``Sphere Packings III'' treats all planar maps in which every face is a
triangle or quadrilateral (except for Boerdijk's pentagonal prism,
treated in ``Sphere Packings V'').
``Sphere Packings IV'' gives preliminary results
on the general planar map.  The main result shows that the planar
map, with a few explicit exceptions, divides the unit sphere
into polygons.  Each polygon is at most a octagon.
The results of
``Sphere Packings IV'' are particularly technical, because
unlike the other cases, we cannot assume that the decomposition stars
have any particular structure.

The preliminaries in IV permit a classification of all planar maps
that are relevant for the proof of the Kepler conjecture.  There
are about 5000 cases that arise.  (This classification was done
by computer.)  After reducing this list to under 100 cases by general
linear programming methods,
the final cases are eliminated case by case.  With
the elimination of the last case, the Kepler conjecture is proved.

\head 8. Computers\endhead

As this project has progressed, the computer has replaced conventional
mathematical arguments more and more, until now
 nearly every aspect of the proof relies on
computer verifications.  Many assertions in these papers
 are results of computer calculations.
To make the proof of Kepler's conjecture more accessible, I have
posted extensive resources \cite{Hal98}.

Computers are used in various significant ways.  They will be
mentioned briefly here, and then developed more thoroughly elsewhere
in the collection, especially in the final paper.

1. {\it  Proof of inequalities by interval arithmetic}.  ``Sphere Packings
I'' describes a method of proving various inequalities in a small number
of variables by computer by interval arithmetic.

2.  {\it Combinatorics}.  A computer program classifies all of the planar maps
that are relevant to the Kepler conjecture.

3. {\it  Linear programming bounds}.  Many of the nonlinear optimization
    problems for the scores of decomposition stars are replaced by linear
    problems that dominate the original score.  They are solved
    by linear programming methods by computer.  A typical problem has
    between 100 and 200 variables and 1000 and 2000 constraints.  Nearly
    100000
    such problems enter into the proof.

4. {\it Branch and bound methods}.  When linear programming methods do not
    give sufficiently good bounds, they have been combined with branch
    and bound methods from global optimization.

5.  {\it Numerical optimization}.  The exploration of the problem
    has been substantially
    aided by nonlinear optimization and symbolic math packages.

6. {\it Organization of output}.
    The organization of the few gigabytes of code and data that
    enter into the proof is in itself a nontrivial undertaking.

\head 9. Acknowledgments\endhead

I am indebted to G. Fejes T\'oth's survey of sphere packings
in the preparation of this overview \cite{Fej97}.
For a much more comprehensive
introduction to the literature on sphere packings, I refer
the reader to that survey and to standard references on
sphere packings such as \cite{CoSl93}, \cite{PaA95}, \cite{Goo97},
\cite{Rog64}, \cite{Fej64}, and \cite{Fej72}.

A detailed strategy of the proof was explained in lectures I gave
at Mount Holyoke and Budapest during the summer of 1996 \cite{Hal96}.
See also the 1996 preprint, ``Recent Progress
on the Kepler Conjecture,'' available from \cite{Hal98}.

I owe the success of this project to a significant degree to
S. Ferguson.  His thesis solves a major step of the program.
He has been highly involved in various other steps of the solution as well.
He returned to Ann Arbor during the final three months of the project
to verify many of the interval-based inequalities appearing in
the appendices of ``Sphere Packings IV'' and ``The Kepler Conjecture.''
It is a pleasure to express my debt to him.

Sean McLaughlin has been involved in this project during the
past year through his fundamental work on the dodecahedral conjecture.
By detecting many of my mistakes, by clarifying my arguments,
 and in many other ways, he has made an important contribution.

 I thank S. Karni, J. Mikhail, J. Song, D. J. Muder, N. J. A. Sloane,
W. Casselman, T. Jarvis, P. Sally, E. Carlson,
and S. Chang for their contributions to this
project.  I express particular thanks to L. Fejes T\'oth for the
inspiration he provided during the course of this research.
This project received the generous institutional support
from
the University of Chicago math department, the
the Insitute for Advanced Study,
the journal {\it Discrete and Computational Geometry},
the School
of Engineering at the University of Michigan (CAEN),
and the National Science
Foundation.
  Software
({\it cfsqp})
\footnote"*"{\quad www.isr.umd.edu/Labs/CACSE/FSQP/fsqp.html}
for testing nonlinear inequalities was provided
by the Institute for Systems Research at the University of Maryland.

Finally, I wish to give my special thanks to Kerri Smith, who has been
my greatest source of support and encouragement through it all.

\newpage
\parindent=0pt

{\bf References}
\smallskip

\parskip=\baselineskip

[Ben74]  Bender, C., Bestimmung der gr\"ossten Anzahl gleich
grosser Kugeln, welche sich auf eine Kugel von demselben
Radius, wie die \"ubrigen, auflegen lassen, {\it Archiv Math.
Physik} 56 (1874), 302--306.

[Bez90] A. Bezdek and W. Kuperberg, Maximum density space packing with
    congruent circular cylinders of infinite length,
    {\it Mathematica} 37 (1990), 74--80.

[BKM91] A. Bezdek, W. Kuperberg, and E. Makai Jr., Maximum density
    space packing with parallel strings of balls,
    {\it DCG} 6 (1991), 227--283.

[Bez94] A. Bezdek, A remark on the packing density in the 3-space
    in {\it Intuitive Geometry}, ed. K. B\"or\"oczky and G. Fejes
    T\'oth, {\it Colloquia Math. Soc. J\'anos Bolyai} 63, North-Holland
    (1994), 17--22.

[Bez97] K. Bezdek, Isoperimetric inequalities and the dodecahedral
    conjecture, {\it Internat. J. Math.} 8, no. 6 (1997), 759--780.

[Bli19] H. F. Blichfeldt,
    Report on the theory of the geometry of numbers,
    {\it Bull. AMS}, 25 (1919), 449--453.

[Bli29] H. F. Blichfeldt,
    The minimum value of quadratic forms and the closest
    packing of spheres, {\it Math. Annalen} 101 (1929), 605--608.

[Bli35] H. F. Blichfeldt,
    The minimum values of positive quadratic forms in six,
    seven and eight variables, {\it Math. Zeit.} 39 (1935), 1--15.

[Boe52] Boerdijk, A. H. Some remarks concerning close-packing
of equal Spheres, {\it Philips Res. Rep.} 7 (1952), 303--313.

[CoHMS94] J. H. Conway, T. C. Hales, D. J. Muder, and N. J. A. Sloane,
    On the Kepler conjecture, {\it Math. Intelligencer} 16,
    no. 2 (1994), 5.

[CoSl95]  J. H. Conway,  N. J. A. Sloane, What are all the
best sphere packings in low dimensions? {\it DCG} 13 (1995), 383--403.

[CoSl93] J. H. Conway and N. J. A. Sloane, Sphere packings, lattices
    and groups,  second edition, Springer-Verlag, New York,
    1993   (third edition, to appear).

[Fej93] G. Fejes T\'oth and W. Kuperberg, Recent results in the
    theory of packing and covering, in New trends in
    discrete and computational geometry, ed. J. Pach, Springer
    1993, 251--279.

[Fej95] G. Fejes T\'oth, Review of [Hsi93], {\it Math. Review} 95g\#52032, 1995.

[Fej95b] G. Fejes T\'oth, Densest packings of typical convex sets
    are not lattice-like, {\it DCG}, 14 (1995), 1--8.

[Fej97] G. Fejes T\'oth, Recent progress on packing and covering,
    preprint.

[Fej72] L. Fejes T\'oth, {\it Lagerungen in der Ebene auf der
    Kugel und im Raum}, second edition,
    Springer-Verlag, Berlin New York, 1972.

[Fej64] L. Fejes T\'oth, Regular figures, Pergamon Press,
    Oxford London New York, 1964.

[Fej42] L. Fejes T\'oth,  \"Uber die dichteste Kugellagerung,
{\it Math. Zeit.} 48 (1942 1943), 676--684.

[Fej50] L. Fejes T\'oth, Some packing and covering theorems,
    {\it Acta Scientiarum Mathematicarum (Szeded)} 12/A, 62--67.

[Fej53] L. Fejes T\'oth, {\it Lagerungen in der Ebene auf
der Kugel und im Raum}, Springer, Berlin, first edition, 1953.

[Gau31] C. F. Gauss, Untersuchungen \"uber die Eigenscahften der
positiven tern\"aren quadratischen Formen von Ludwig August Seber,
    {\it G\"ottingische gelehrte Anzeigen}, 1831 Juli 9,
also published in {\it J. reine angew. Math.} 20 (1840), 312--320, and
    {\it Werke},  vol. 2,
    K\"onigliche Gesellschaft der Wissenschaften, G\"ottingen,
            1876, 188--196.

[Goo97] J. E. Goodman and J. O'Rourke, Handbook of discrete and
    computational geometry, CRC, Boca Raton and New York, 1997.

[Gun75] S. G\"unther, {\it Ein stereometrisches Problem},
{\it Archiv der Math. Physik} 57 (1875), 209--215.

[Hal92] T. C. Hales, The sphere packing problem, {\it J. Computational
    Applied Math.} 44 (1992), 41--76.

[Hal93] T. C. Hales, Remarks on the density of sphere packings in
    three dimensions, {\it Combinatorica} 13 (1993), 181--187.

[Hal94] T. C. Hales, The status of the Kepler conjecture,
    {\it Math. Intelligencer} 16, no. 3, (1994), 47--58.

[Hal96] T. C. Hales,
    {http://www.math.lsa.umich.edu/\~\relax hales/holyoke.html}

[Hal98] T. C. Hales, \hfill{http://www.math.lsa.umich.edu/\~%
    \relax hales/packings.html}

[Hil01] D. Hilbert, Mathematische Probleme, {\it Archiv Math. Physik} 1 (1901),
    44--63, also in {\it Proc. Sym. Pure Math.} 28 (1976), 1--34.

[Hop74] Hoppe R. {\it Bemerkung der Redaction}, Math. Physik 56
(1874), 307-312.

[HoPT95] R. Horst, P.M. Pardalos, N.V. Thoai, {\it Introduction
    to Global Optimization}, Kluwer, 1995.

[Hsi93] W.-Y. Hsiang, On the sphere packing problem and the proof
    of Kepler's conjecture, Internat. J. Math 93 (1993), 739-831.

[Hsi93a] W.-Y. Hsiang, On the sphere packing problem and the
    proof of Kepler's conjecture, in {\it Differential geometry and
    topology} (Alghero, 1992), World Scientific, River Edge,
    NJ, 1993,  117--127.

[Hsi93b] W.-Y. Hsiang, The geometry of spheres, in {\it Differential
    geometry} (Shanghai, 1991), World Scientific, River Edge, NJ,
    1993, 92-107.

[Hsi95] W.-Y. Hsiang, A rejoinder to T. C. Hales's article ``The status
    of the Kepler conjecture,'' {\it Math. Intelligencer} 17, no. 1, (1995),
    35--42.

[Kar66] R. Kargon, Atomism in England from Hariot to Newton,
    Oxford, 1966.

[Kep66] J. Kepler, The Six-cornered snowflake, Oxford Clarendon Press,
    Oxford, 1966,  forward by L. L. Whyte.

[KoZ73] A. Korkine and  G. Zolotareff, Sur les formes quadratiques,
    {\it Math. Annalen} 6 (1873), 366--389.

[KoZ77] A. Korkine and  G. Zolotareff, Sur les formes quadratiques
    positives, {\it Math. Annalen} 11 (1877), 242--292.

[Lag73] J. L. Lagrange,  Recherches d'arithm\'etique, {\it Nov. Mem.
    Acad. Roy. Sc. Bell Lettres Berlin} 1773, in {\it \OE uvres}, vol. 3,
    693--758.

[Lee56] J. Leech, The Problem of the Thirteen Spheres,
{\it The Mathematical Gazette}, Feb 1956, 22--23.

[Lin86] J. H. Lindsey II, Sphere packing in $R^3$, {\it Mathematika}
    33 (1986), 137--147.

[Mas66] B. J. Mason, On the shapes of snow crystals, in \cite{Kep66}.

[McL98] S. McLaughlin, A proof of the dodecahedral conjecture,
    preprint.

[Mel97] J. B. M. Melissen, Packing and covering with circles,
    Ph.D. dissertation, Univ. Utrecht, Dec. 1997.

[Mil76] J. Milnor, Hilbert's problem 18: on crystallographic groups,
    fundamental domains, and on sphere packings, in
    Mathematical developments arising from Hilbert problems,
    {\it Proc. Symp. Pure Math.}, vol 28, 491--506, AMS, 1976.

[MP93] W. Moser, J. Pach, Research problems in discrete geometry,
    DIMACS Technical Report, 93032, 1993.

[Mu88] D. J. Muder, Putting the best face on a Voronoi polyhedron,
    {\it Proc. London Math. Soc.} (3) 56 (1988), 329--348.

[Mud93]  D. J. Muder A New Bound on the Local Density
of Sphere Packings, {\it Discrete and Comp. Geom.} 10 (1993), 351--375.

[Mud97]  D. J. Muder, letter, in {\it Fermat's enigma}, by S. Singh,
        Walker, New York, 1997.

[Oes90] J. Oesterl\'e,  Empilements de sph\`eres,
    S\'eminaire Bourbaki, vol. 1989/90, Ast\'erisque (1990),
        No. 189--190 exp. no. 727, 375--397.

[PaA95] J. Pach, P.K. Agarwal, {\it Combinatorial geometry}, John Wiley,
    New York 1995.

[Plo00]  K. Plofker, private communication, January 2000.

[Ran47] R. A. Rankin, {\it Annals of Math.} 48 (1947), 228--229.

[Rog58] C. A. Rogers, The packing of equal spheres, {\it Proc. London Math.
    Soc.} (3) 8 (1958), 609--620.

[Rog64] C. A. Rogers, {\it Packing and covering}, Cambridge University Press,
    Cambridge, 1964.

[Sch53] K. Sch\"utte and B.L. van der Waerden, Das
Problem der dreizehn Kugeln, {\it Math. Annalen} 125, (1953), 325--334.

[SeM44] B. Segre and K. Mahler, On the densest packing of
    circles, {\it Amer. Math Monthly} (1944), 261--270.

[Shi83] J. W. Shirley,
{\it Thomas Harriot: a biography}, Oxford, 1983.

[SHDC95] N. J. A. Sloane, R. H. Hardin, T. D. S. Duff, J. H. Conway,
    Minimal-energy clusters of hard spheres,
    {\it DCG} 14,  no. 3, (1995), 237--259.

[Thu92] A. Thue, Om nogle geometrisk taltheoretiske Theoremer,
    {\it Forandlingerneved de Skandinaviske Naturforskeres} 14 (1892), 352--353.

[Thu10] A. Thue, \"Uber die dichteste Zusammenstellung von
    kongruenten Kreisen in der Ebene, {\it Christinia Vid. Selsk. Skr.} 1
    (1910), 1--9.

[Why66] L. L. Whyte, forward to \cite{Kep66}.

\bye